\documentclass[12pt]{article}
\usepackage{amssymb}

\renewcommand{\labelenumi}{(\theenumi)}
\renewcommand{\theenumi}{\roman{enumi}}
\newtheorem{definition}{Definition}[section]
\newtheorem{theorem}[definition]{Theorem}
\newtheorem{lemma}[definition]{Lemma}
\newtheorem{corollary}[definition]{Corollary}
\newtheorem{remark}[definition]{Remark}
\newtheorem{example}[definition]{Example}

\newtheorem{note}[definition]{Note}

\def\R{\mathbb R}
\def\C{\mathbb C}
\def\K{\mathbb K}
\def\D{\mathbb D}
\def\Z{\mathbb Z}
\def\K{\mathbb K}
\def\fld{\mathbb K}
\def\Mdf{{\mbox{Mat}}_{d+1}(\K)}

\begin{document}

\title{ \bf The equitable presentation
for the \\ quantum  
group $U_q({\mathfrak{g}})$ 
associated with
a  \\symmetrizable Kac-Moody algebra 
$\mathfrak g$
}
\author{
Paul Terwilliger{\footnote{
Department of Mathematics, University of
Wisconsin, 480 Lincoln Drive, Madison WI 53706-1388 USA}
}}
\date{}
%to get date printout, comment out above line 

\maketitle
\begin{abstract} 
We consider the quantum group
$U_q(\mathfrak{g})$ associated
with a symmetrizable Kac-Moody algebra 
${\mathfrak{g}}$.
We display a presentation for 
$U_q(\mathfrak{g})$ that
we find attractive; we call this 
the 
{\it equitable}
presentation.
 For  
$\mathfrak{g}=
\mathfrak{sl}_2$ the
equitable presentation
has generators $X^{\pm 1},Y,Z$
and relations $X X^{-1} = X^{-1} X=1$,
\begin{eqnarray*}
\frac{qXY-q^{-1}YX}{q-q^{-1}}=1,
\qquad 
\frac{qYZ-q^{-1}ZY}{q-q^{-1}}=1,
\qquad 
\frac{qZX-q^{-1}XZ}{q-q^{-1}}=1.
\end{eqnarray*}
\end{abstract}

\section{Introduction}
In \cite{equit} T. Ito, C. Weng, and the present
author
introduced the equitable presentation
for the  quantum group 
$U_q(\mathfrak{sl}_2)$. 
The purpose of this
note is to give an analogous presentation
for 
the quantum group
$U_q(\mathfrak{g})$, where
$\mathfrak{g}$ is a symmetrizable
Kac-Moody algebra. 
As we will see,
the generators for this presentation
are related to
Koornwinder's 
twisted primitive elements 
\cite{koorn1},
\cite{koorn}.

\medskip
\noindent 
Throughout, $\K$ will denote a field and
 $q$ will denote an indeterminate. 
We will work over the field $\K(q)$.

\medskip
\noindent In order to motivate our
main result we first recall some facts
about
$U_q(\mathfrak{sl}_2)$.

\section{The quantum group
$U_q(\mathfrak{sl}_2)$}

\noindent We begin with the definition of 
$U_q(\mathfrak{sl}_2)$.

\begin{definition} 
\cite[p.~9]{Jantzen}
\label{def:uq}
\rm
We let $U_q(\mathfrak{sl}_2)$
denote the unital associative $\K(q)$-algebra 
with
generators $E,F,K, K^{-1}$
and the following relations:
\begin{eqnarray*}
KK^{-1} &=& 
K^{-1}K =  1,
\label{eq:buq1}
\\
KEK^{-1} &=& q^2E,
\label{eq:buq2}
\\
KFK^{-1}&=&q^{-2}F,
\label{eq:buq3}
\\
EF-FE &=& \frac{K-K^{-1}}{q-q^{-1}}.
\label{eq:buq4}
\end{eqnarray*}
\end{definition}

\noindent 
We call $E,F,K,K^{-1}$
the {\it Chevalley generators} for
$U_q(\mathfrak{sl}_2)$.

\medskip
\noindent 
We refer the reader to
\cite{cp},
\cite{Jantzen},
\cite{Kassel}
 for
background information on 
$U_q(\mathfrak{sl}_2)$.

\begin{lemma} 
 \cite[p.~35]{Jantzen}
\label{lem:hopfsl2}
The quantum group
$U_q(\mathfrak{sl}_2)$
has the following Hopf algebra structure. 
The comultiplication $\Delta$ satisfies
\begin{eqnarray*}
\Delta(E) &=& E \otimes 1+K \otimes E,
\\
\Delta(F) &=& F \otimes K^{-1}+ 1 \otimes F,
\\
\Delta(K) &=& K \otimes K.
\end{eqnarray*}
The counit $\varepsilon$ satisfies
\begin{eqnarray*}
\varepsilon(E)=0, \qquad \qquad \varepsilon(F)=0,
\qquad \qquad 
\varepsilon(K)=1.
\end{eqnarray*}
The antipode $S$ satisfies
\begin{eqnarray*}
S(E)= -K^{-1}E,
\qquad \qquad 
S(F)= -FK,
\qquad \qquad 
S(K) = K^{-1}.
\end{eqnarray*}
\end{lemma}

\section{The equitable presentation for
$U_q(\mathfrak{sl}_2)$
}

In
the presentation for 
$U_q(\mathfrak{sl}_2)$ given
in 
Definition
\ref{def:uq}
the generators $E,F$
and the generators $K, K^{-1}$ play a very different
role. 
We now recall a presentation for 
$U_q(\mathfrak{sl}_2)$ whose
generators
are on a more equal footing.  

\begin{theorem}
\cite[Theorem 2.1]{equit}
\label{thm:uq2}
The algebra
$U_q(\mathfrak{sl}_2)$ 
 is isomorphic to
the unital associative $\K(q)$-algebra 
with
generators 
$X$, $X^{-1}$, $Y$, $Z$
and the following relations:
\begin{eqnarray*}
XX^{-1} = 
X^{-1}X &=&  1,
\label{eq:2buq1}
\\
\frac{qXY-q^{-1}YX}{q-q^{-1}}&=&1,
\label{eq:2buq2}
\\
\frac{qYZ-q^{-1}ZY}{q-q^{-1}}&=&1,
\label{eq:2buq3}
\\
\frac{qZX-q^{-1}XZ}{q-q^{-1}}&=&1.
\label{eq:2buq4}
\end{eqnarray*}
An isomorphism with the presentation in Definition
\ref{def:uq} is given by:
\begin{eqnarray*}
\label{eq:iso1}
X^{{\pm}1} &\rightarrow & K^{{\pm}1},\\
\label{eq:iso2}
Y &\rightarrow & K^{-1}+F(q-q^{-1}), \\
\label{eq:iso3}
Z &\rightarrow & K^{-1}-K^{-1}Eq(q-q^{-1}).
\end{eqnarray*}
The inverse of this isomorphism is given by:
\begin{eqnarray*}
\label{eq:iso3inv}
E &\rightarrow & (1-XZ)q^{-1}(q-q^{-1})^{-1},\\
\label{eq:iso2inv}
F &\rightarrow & (Y-X^{-1})(q-q^{-1})^{-1}, \\
\label{eq:iso1inv}
K^{{\pm}1} &\rightarrow & X^{{\pm}1}.\\
\end{eqnarray*}
\end{theorem}
\noindent {\it Proof:} One readily checks that
each map is a homomorphism of $\K(q)$-algebras and that
the maps are inverses.
It follows that each map is an isomorphism of $\K(q)$-algebras.
\hfill $\Box $ \\

\noindent
The generators 
 $X$, $Y$, $Z$ from
Theorem
\ref{thm:uq2} 
 are 
on an equal footing, more or less.
This motivates the following definition.

\begin{definition}
\label{def:equit2}
\cite[Definition 2.2]{equit}
\rm
By the {\it equitable} presentation for
$U_q(\mathfrak{sl}_2)$ we mean
the presentation given in
 Theorem
\ref{thm:uq2}.
We call 
$X,X^{-1}, Y,Z$ the
{\it equitable generators}.
\end{definition}

\begin{definition}
\label{def:nc}
\rm
For notational convenience, throughout this paper
we identify
the copy of 
$U_q(\mathfrak{sl}_2)$ given 
in Definition
\ref{def:uq} with the copy of
$U_q(\mathfrak{sl}_2)$ given 
in Theorem
\ref{thm:uq2}, via the isomorphism
given 
in Theorem
\ref{thm:uq2}.
\end{definition}

\noindent The Hopf algebra
structure for 
$U_q(\mathfrak{sl}_2)$ given in 
Lemma
\ref{lem:hopfsl2} looks as follows
in terms of the 
equitable generators.

\begin{theorem}
With reference to
Lemma
\ref{lem:hopfsl2}  and
Definition \ref{def:nc},
the comultiplication $\Delta$ 
 satisfies
\begin{eqnarray*}
\Delta(X) &=& X \otimes X,
\\
\Delta(Y) &=& 
(Y-1) \otimes X^{-1}+
1 \otimes Y,
\\
\Delta(Z) &=& 
(Z-1) \otimes X^{-1}+
1 \otimes Z.
\end{eqnarray*}
The counit $\varepsilon$ 
satisfies
\begin{eqnarray*}
\varepsilon(X)=1, \qquad \qquad \varepsilon(Y)=1,
\qquad \qquad 
\varepsilon(Z)=1.
\end{eqnarray*}
The antipode $S$ 
satisfies
\begin{eqnarray*}
S(X)= X^{-1},
\qquad \qquad 
S(Y)= 1+X-YX,
\qquad \qquad 
S(Z) = 1+X-ZX.
\end{eqnarray*}
\end{theorem}
\noindent {\it Proof:} 
Routine verification.
\hfill $\Box $ \\

\noindent We finish this section with a remark.

\begin{remark}
\label{koorn}
\rm
With reference
to 
Lemma
\ref{lem:hopfsl2} and
Definition \ref{def:nc},
for $y=Y-1$ and $z=Z-1$ we have
\begin{eqnarray*}
\Delta(y) &=& 
y \otimes X^{-1}+
1 \otimes y,
\\
\Delta(z) &=& 
z\otimes X^{-1}+
1 \otimes z.
\end{eqnarray*}
We will discuss these two equations in
Section 6.
\end{remark}

\section{The quantum group
$U_q(\mathfrak{g})$
}

\noindent 
We now turn our attention to the quantum
group 
$U_q(\mathfrak{g})$, where
$\mathfrak{g}$ is a symmetrizable Kac-Moody
algebra.
We will give an ``equitable'' presentation
for 
$U_q(\mathfrak{g})$.
We begin with a comment. In the standard references
\cite{brown},
\cite{cp},
\cite{HK},
\cite{Jantzen},
\cite{joseph},
\cite{Kassel},
\cite{ks},
\cite{lustig},
\cite{majid}
the definition 
of 
$U_q(\mathfrak{g})$ varies. We will
give the equitable presentation for
one version
of 
$U_q(\mathfrak{g})$; the equitable presentation
for the other versions can be 
obtained from this one
by minor 
modification.

\medskip
\noindent 
We will use the following notation.

\begin{definition}
\label{def:prelim}
\rm
Let $n$ denote a positive integer and
let $A$ denote a symmetrizable generalized Cartan 
matrix\cite[p.~1]{kac} 
of order $n$.
Since $A$ is symmetrizable
%by 
%\cite[p.~16]{kac}
 there exists relatively prime positive integers
$d_1,\ldots, d_n$ such that $d_i A_{ij}=d_jA_{ji}$
for $1 \leq i,j\leq n$.
For $1 \leq i \leq n$ we do the following.
 We define 
$q_i = q^{d_i}$.
Also, for an integer  $m$ we define 
\begin{eqnarray*}
\lbrack m \rbrack_i = \frac {q_i^m-q_i^{-m}}{q_i-q_i^{-1}}
\end{eqnarray*}
and for $m\geq 0$ we define
\begin{eqnarray*}
\lbrack m \rbrack_i^{!} = 
\lbrack m \rbrack_i  \lbrack m-1 \rbrack_i
\cdots \lbrack 2 \rbrack_i \lbrack 1 \rbrack_i.
\end{eqnarray*}
We interpret $\lbrack 0 \rbrack_i^!=1$.
For integers $m\geq  r \geq 0$
we define
\begin{eqnarray*}
\biggl\lbrack {{m}\atop{r}} \biggr\rbrack_i
= \frac{\lbrack m \rbrack_i^!}{\lbrack r \rbrack_i^! \lbrack m-r \rbrack_i^!}.
\end{eqnarray*}
Let
$\mathfrak{g}=
\mathfrak{g}'(A)$
 denote the 
Kac-Moody algebra over $\C$ that corresponds to $A$
\cite[p.~xi]{kac}. 
\end{definition}

\begin{definition} 
\label{def:uqg}
\rm
\cite[p.~281]{cp}
%\cite[p.~262]{cp}
With reference  to Definition \ref{def:prelim},
$U_q(\mathfrak{g})$  is the 
unital associative $\K(q)$-algebra 
with
generators 
$E_i, F_i,
K_i, K^{-1}_i$
$(i=1,\ldots,n)$
and the following relations:
\begin{enumerate}
\item[{\rm (R1)}]
$K_iK^{-1}_i = 
K^{-1}_iK_i =  1$
\item[{\rm (R2)}]
$K_iK_j = K_jK_i$
\item[{\rm (R3)}]
$K_iE_jK^{-1}_i=q_i^{A_{ij}}E_j$
\item[{\rm (R4)}]
$K_iF_jK^{-1}_i=q_i^{-A_{ij}}F_j$
\item[{\rm (R5)}]
$\displaystyle{
E_iF_j-F_jE_i = \delta_{ij}\frac{K_i-K_i^{-1}}{q_i-q^{-1}_i}}$
\item[{\rm (R6)}]
$\displaystyle{
\sum_{r=0}^{1-A_{ij}} (-1)^r 
\biggl\lbrack {{1-A_{ij}}\atop{r}} \biggr\rbrack_i
E^{1-A_{ij}-r}_i E_j E^r_i =0 \qquad \mbox{if} \qquad i\not=j}$
\item[{\rm (R7)}]
$\displaystyle{
\sum_{r=0}^{1-A_{ij}} (-1)^r 
\biggl\lbrack {{1-A_{ij}}\atop{r}} \biggr\rbrack_i
F^{1-A_{ij}-r}_i F_j F^r_i =0 \qquad \mbox{if} \qquad i\not=j}$.
\end{enumerate}
The expression
 $\delta_{ij}$
in (R5) is the Kronecker delta.
\end{definition}

\noindent 
We call 
the generators $E_i, F_i, K_i, K^{-1}_i$ from
Definition
\ref{def:uqg}
the {\it Chevalley generators} for
$U_q(\mathfrak{g})$.

%\medskip
%\noindent 
%We refer the reader to
% \cite{}
% for
%background information on 
%$U_q({\widehat{\mathfrak{g}}})$.

\medskip
\noindent 
We now recall a Hopf algebra structure for 
$U_q(\mathfrak{g})$.

\begin{lemma} \cite[p.~55]{Jantzen}
\label{lem:hopfuqg}
The quantum group
$U_q(\mathfrak{g})$
has the following Hopf algebra structure. 
The
 comultiplication $\Delta$ satisfies
\begin{eqnarray*}
\Delta(E_i) &=& E_i \otimes 1+K_i \otimes E_i,
\\
\Delta(F_i) &=& F_i \otimes K_i^{-1}+ 1 \otimes F_i,
\\
\Delta(K_i) &=& K_i \otimes K_i.
\end{eqnarray*}
The counit $\varepsilon$ satisfies
\begin{eqnarray*}
\varepsilon(E_i)=0, \qquad \qquad \varepsilon(F_i)=0,
\qquad \qquad 
\varepsilon(K_i)=1.
\end{eqnarray*}
The antipode $S$ satisfies
\begin{eqnarray*}
S(E_i)= -K_i^{-1}E_i,
\qquad \qquad 
S(F_i)= -F_iK_i,
\qquad \qquad 
S(K_i) = K_i^{-1}.
\end{eqnarray*}
\end{lemma}

\section{The equitable presentation for
$U_q(\mathfrak{g})$}

\noindent In this section we give the
equitable presentation for
$U_q(\mathfrak{g})$. The following is our main
result.

\begin{theorem}
\label{thm:uqg2}
The
$\K(q)$-algebra
$U_q(\mathfrak{g})$
from
Definition
\ref{def:uqg}
is isomorphic to
the unital associative $\K(q)$-algebra 
with
generators
$X_i, X^{-1}_i, Y_i, Z_i$
$(i=1,\ldots, n)$
and the following relations:
\begin{enumerate}
\item[{\rm (E1)}]
$X_iX^{-1}_i = 
X^{-1}_iX_i =  1$
\item[{\rm (E2)}]
$X_iX_j = X_jX_i$
\item[{\rm (E3)}]
$Y_iX_j-q_i^{A_{ij}}X_jY_i=
X^{-1}_iX_j(1-q_i^{A_{ij}})$
\item[{\rm (E4)}]
$X_iZ_j-q_i^{A_{ij}}Z_jX_i=
X_iX^{-1}_j(1-q_i^{A_{ij}})$
\item[{\rm (E5)}]
$Z_iY_i-q_i^2 Y_iZ_i=
1-q_i^2$
\item[{\rm (E6)}]
$Z_iY_j-q_i^{A_{ij}}Y_jZ_i=
X^{-1}_iX^{-1}_j(1-q_i^{A_{ij}}) \qquad \mbox{if} \quad i \not=j$
\item[{\rm (E7)}]
$\displaystyle{
\sum_{r=0}^{1-A_{ij}} (-1)^r 
\biggl\lbrack {{1-A_{ij}}\atop{r}} \biggr\rbrack_i
Y^{1-A_{ij}-r}_i Y_j Y^r_i =
X_i^{A_{ij}-1}X^{-1}_j \prod_{s=0}^{-A_{ij}}(1-q_i^{A_{ij}+2s})
\qquad \mbox{if} \qquad i\not=j}$
\item[{\rm (E8)}]
$\displaystyle{
\sum_{r=0}^{1-A_{ij}} (-1)^r 
\biggl\lbrack {{1-A_{ij}}\atop{r}} \biggr\rbrack_i
Z^{1-A_{ij}-r}_i Z_j Z^r_i =
X_i^{A_{ij}-1}X^{-1}_j \prod_{s=0}^{-A_{ij}}(1-q_i^{A_{ij}+2s})
\qquad \mbox{if} \qquad i\not=j}$
\end{enumerate}
An isomorphism with the presentation for
$U_q(\mathfrak{g})$ given in
Definition
\ref{def:uqg} is
\begin{eqnarray}
\label{eq:i1}
X^{\pm 1}_i &\rightarrow & K^{\pm 1}_i,
\\
\label{eq:i2}
Y_i&\rightarrow & K^{-1}_i +F_i(q_i-q_i^{-1}),
\\
\label{eq:i3}
Z_i&\rightarrow &K^{-1}_i -K_i^{-1}E_iq_i(q_i-q_i^{-1}).
\end{eqnarray}
The inverse of this isomorphism is
\begin{eqnarray}
\label{eq:iinv1}
E_i&\rightarrow &(1-X_iZ_i)q^{-1}_i(q_i-q_i^{-1})^{-1},
\\
\label{eq:iinv2}
F_i&\rightarrow &(Y_i-X^{-1}_i)(q_i-q_i^{-1})^{-1},
\\
\label{eq:iinv3}
K^{\pm 1}_i &\rightarrow & X^{\pm 1}_i.
\end{eqnarray}
\end{theorem}
\noindent {\it Proof:}
We recall a few identities.
By \cite[p.~6]{Jantzen} we find that
for integers $m\geq r\geq 1$,
\begin{eqnarray}
\label{eq:3ident}
\biggl\lbrack {{m}\atop{r}} \biggr\rbrack_i
+q^{m+1}_i\biggl\lbrack {{m}\atop{r-1}} \biggr\rbrack_i
=q_i^r\biggl\lbrack {{m+1}\atop{r}} \biggr\rbrack_i
\qquad \qquad (1 \leq i \leq n).
\end{eqnarray}
The following 
identity is a special case of the
$q$-binomial theorem \cite[p.~236]{GR}.
For an indeterminate $\lambda$
and for an integer $m\geq 0$,
\begin{eqnarray}
\label{eq:vand2}
\sum_{r=0}^m (-1)^r
\biggl\lbrack {{m}\atop{r}} \biggr\rbrack_i
\lambda^r &=& 
 \prod_{s=0}^{m-1}(1-\lambda q_i^{1-m+2s})
\qquad \qquad (1 \leq i \leq n).
\end{eqnarray}
Let 
$U$ denote
the unital associative $\K(q)$-algebra 
with
generators
$X_i, X^{-1}_i, Y_i, Z_i$
$(i=1,\ldots, n)$
and relations (E1)--(E8).
Our goal is to show
that (\ref{eq:i1})--(\ref{eq:i3}) gives
an isomorphism
of 
$\K(q)$-algebras from
$U$ to
$U_q(\mathfrak{g})$, and that the
inverse of this isomorphism satisfies 
(\ref{eq:iinv1})--(\ref{eq:iinv3}).
To this end, we first show that 
(\ref{eq:i1})--(\ref{eq:i3}) gives
a homomorphism of 
$\K(q)$-algebras
from
$U$ to
$U_q(\mathfrak{g})$.
In order to do this we let
${\hat X}_i^{\pm 1},
{\hat Y}_i, {\hat Z}_i$ denote the expressions
on the right in 
(\ref{eq:i1})--(\ref{eq:i3}) respectively
and show the following hold in
$U_q(\mathfrak{g})$:
\begin{eqnarray}
&&{\hat X}_i{\hat X}^{-1}_i = 
{\hat X}^{-1}_i{\hat X}_i =  1,
\label{eq:hat1}
\\
&&
{\hat X}_i{\hat X}_j = {\hat X}_j{\hat X}_i,
\label{eq:hat2}
\\
&&
{\hat Y}_i{\hat X}_j-q_i^{A_{ij}}{\hat X}_j{\hat Y}_i=
{\hat X}^{-1}_i{\hat X}_j(1-q_i^{A_{ij}}),
\label{eq:hat3}
\\
&&{\hat X}_i{\hat Z}_j-q_i^{A_{ij}}{\hat Z}_j{\hat X}_i=
{\hat X}_i{\hat X}^{-1}_j(1-q_i^{A_{ij}}),
\label{eq:hat4}
\\
&&
{\hat Z}_i{\hat Y}_i-q_i^2 {\hat Y}_i{\hat Z}_i=
1-q_i^2,
\label{eq:hat5}
\\
&&
{\hat Z}_i{\hat Y}_j-q_i^{A_{ij}}{\hat Y}_j{\hat Z}_i=
{\hat X}^{-1}_i{\hat X}^{-1}_j(1-q_i^{A_{ij}})
\quad \mbox{if} \quad i \not=j,
\label{eq:hat6}
\\
&&
\sum_{r=0}^{1-A_{ij}} (-1)^r 
\biggl\lbrack {{1-A_{ij}}\atop{r}} \biggr\rbrack_i
{\hat Y}^{1-A_{ij}-r}_i {\hat Y}_j {\hat Y}^r_i =
{\hat X}_i^{A_{ij}-1}{\hat X}^{-1}_j
\prod_{s=0}^{-A_{ij}}(1-q_i^{A_{ij}+2s})
\quad \mbox{if} \quad i\not=j
\qquad  
\label{eq:hat7}
\\
&&
\sum_{r=0}^{1-A_{ij}} (-1)^r 
\biggl\lbrack {{1-A_{ij}}\atop{r}} \biggr\rbrack_i
{\hat Z}^{1-A_{ij}-r}_i {\hat Z}_j {\hat Z}^r_i =
{\hat X}_i^{A_{ij}-1}{\hat X}^{-1}_j 
\prod_{s=0}^{-A_{ij}}(1-q_i^{A_{ij}+2s})
\quad \mbox{if} \quad i\not=j.
\label{eq:hat8}
\end{eqnarray}
Line
(\ref{eq:hat1}) holds by 
(R1) and line
(\ref{eq:hat2}) holds by 
(R2).
Line 
(\ref{eq:hat3}) follows
from
(R1), (R2), (R4).
Line 
(\ref{eq:hat4}) follows
from
(R1), (R2), (R3).
Lines 
(\ref{eq:hat5}), (\ref{eq:hat6}) each
follow
from
(R1)--(R5).
In order to verify
(\ref{eq:hat7}) we use the following identity.
By (R4) and induction we obtain
\begin{eqnarray}
\label{eq:b1}
(K^{-1}_i + F_i(q_i-q^{-1}_i))^m
=
\sum_{r=0}^m 
\biggl\lbrack {{m}\atop{r}} \biggr\rbrack_i
q_i^{r(m-r)}(q_i-q^{-1}_i)^rF_i^rK_i^{r-m}
\end{eqnarray}
for all integers $m\geq 0$.
In the equation of (\ref{eq:hat7})
the expression on the left is equal to
\begin{eqnarray}
\label{eq:b2}
\sum_{r=0}^{1-A_{ij}} (-1)^r 
\biggl\lbrack {{1-A_{ij}}\atop{r}} \biggr\rbrack_i
(K_i^{-1}+F_i(q_i-q_i^{-1}))^{1-A_{ij}-r} 
K_j^{-1}
(K_i^{-1}+F_i(q_i-q_i^{-1}))^r
\end{eqnarray}
plus $q_j-q_j^{-1}$ times 
\begin{eqnarray}
\label{eq:b3}
\sum_{r=0}^{1-A_{ij}} (-1)^r 
\biggl\lbrack {{1-A_{ij}}\atop{r}} \biggr\rbrack_i
(K_i^{-1}+F_i(q_i-q_i^{-1}))^{1-A_{ij}-r} 
F_j
(K_i^{-1}+F_i(q_i-q_i^{-1}))^r,
\end{eqnarray}
and the expression on the right is
\begin{eqnarray}
K_i^{A_{ij}-1}K_j^{-1}
\prod_{s=0}^{-A_{ij}}(1-q_i^{A_{ij}+2s}).
\label{eq:b4}
\end{eqnarray}
For $i \not=j$ we show the expressions
(\ref{eq:b2}),  
(\ref{eq:b4}) are equal and 
the expression (\ref{eq:b3}) is zero.
In 
(\ref{eq:b2}) we evaluate
the terms
$(K_i^{-1}+F_i(q_i-q_i^{-1}))^{1-A_{ij}-r}$
and
$(K_i^{-1}+F_i(q_i-q_i^{-1}))^r
$
using
(\ref{eq:b1}); the result is
\begin{eqnarray*}
&&\sum_{r=0}^{1-A_{ij}} 
\sum_{\eta=0}^{1-A_{ij}-r} 
\sum_{t=0}^r 
(-1)^r 
\biggl\lbrack {{1-A_{ij}}\atop{r}} \biggr\rbrack_i
\biggl\lbrack {{1-A_{ij}-r}\atop{\eta}} \biggr\rbrack_i
\biggl\lbrack {r\atop{t}} \biggr\rbrack_i
q_i^{\eta(1-A_{ij}-\eta-r)}
\nonumber
\\
\label{eq:b5}
&&
\qquad \qquad \qquad \times \quad 
q_i^{t(r-t)}
(q_i-q_i^{-1})^{\eta+t}
F_i^\eta 
K_i^{A_{ij}+\eta+r-1}
K_j^{-1}
F_i^t 
K_i^{t-r}.
\end{eqnarray*}
In 
the above expression
we 
pull the
terms
$K_i^{A_{ij}+\eta+r-1}
K_j^{-1}
$
to the right past the $F_i$
using (R4). We obtain

\parbox{10cm}
{\begin{eqnarray*}
&&\sum_{r=0}^{1-A_{ij}} 
\sum_{\eta=0}^{1-A_{ij}-r} 
\sum_{t=0}^r 
(-1)^r 
\biggl\lbrack {{1-A_{ij}}\atop{r}} \biggr\rbrack_i
\biggl\lbrack {{1-A_{ij}-r}\atop{\eta}} \biggr\rbrack_i
\biggl\lbrack {r\atop{t}} \biggr\rbrack_i
q_i^{(\eta+t)(1-A_{ij}-\eta-r)}
%\nonumber
\\
%\label{eq:b6}
&&
\qquad \qquad \qquad \qquad \qquad \times \quad 
q_i^{t(1-\eta-t)}
(q_i-q_i^{-1})^{\eta+t}
F_i^{\eta+t}K_i^{A_{ij}+\eta+t-1}K_j^{-1}.
\end{eqnarray*}}
\hfill
\parbox{1cm}{\begin{eqnarray} 
\label{eq:b6}
\end{eqnarray}}

\noindent In 
(\ref{eq:b6})
we change variables by substituting
$u=r-t$ and $v=\eta+t$. We find
that for $0 \leq v \leq 1-A_{ij}$
the coefficient of
$F_i^vK_i^{A_{ij}+v-1}K_j^{-1}$
in 
(\ref{eq:b6}) is equal to
\begin{eqnarray}
\biggl\lbrack {{1-A_{ij}}\atop{v}} \biggr\rbrack_i
q_i^{v(1-A_{ij}-v)}
(q_i-q^{-1}_i)^v
\label{eq:b7}
\end{eqnarray}
times
\begin{eqnarray}
\sum_{t=0}^v 
(-1)^t 
\biggl\lbrack {{v}\atop{t}} \biggr\rbrack_i
q_i^{t(1-v)}
\label{eq:b8}
\end{eqnarray}
times
\begin{eqnarray}
\sum_{u=0}^{1-A_{ij}-v} 
(-1)^{u} 
\biggl\lbrack {{1-A_{ij}-v}\atop{u}} \biggr\rbrack_i
q_i^{-uv}.
\label{eq:b9}
\end{eqnarray}
The sum
(\ref{eq:b8}) is equal to 
$\delta_{v,0}$
in view of
 (\ref{eq:vand2}).
For $v=0$
the expression
(\ref{eq:b7}) is equal to 1
and the sum 
(\ref{eq:b9}) is equal to
\begin{eqnarray}
\prod_{s=0}^{-A_{ij}}(1-q_i^{A_{ij}+2s})
\label{eq:b10}
\end{eqnarray}
in view of
(\ref{eq:vand2}).
Therefore 
(\ref{eq:b6}) is equal to
(\ref{eq:b4}).
From our above comments
(\ref{eq:b2}) is equal to
(\ref{eq:b4}).
In 
(\ref{eq:b3}) we evaluate
the terms
$(K_i^{-1}+F_i(q_i-q_i^{-1}))^{1-A_{ij}-r}$
and
$(K_i^{-1}+F_i(q_i-q_i^{-1}))^r
$
using
(\ref{eq:b1}); the result is
\begin{eqnarray*}
&&\sum_{r=0}^{1-A_{ij}} 
\sum_{\eta=0}^{1-A_{ij}-r} 
\sum_{t=0}^r 
(-1)^r 
\biggl\lbrack {{1-A_{ij}}\atop{r}} \biggr\rbrack_i
\biggl\lbrack {{1-A_{ij}-r}\atop{\eta}} \biggr\rbrack_i
\biggl\lbrack {r\atop{t}} \biggr\rbrack_i
q_i^{\eta(1-A_{ij}-\eta-r)}
\\
&&
\qquad \qquad \qquad \times \quad 
q_i^{t(r-t)}
(q_i-q_i^{-1})^{\eta+t}
F_i^\eta 
K_i^{A_{ij}+\eta+r-1}
F_j
F_i^t 
K_i^{t-r}.
\end{eqnarray*}
\noindent In the above expression
we 
pull the
terms
$K_i^{A_{ij}+\eta+r-1}
$
to the right past the $F_j$, $F_i$
using (R4). We obtain

\parbox{10cm}
{\begin{eqnarray*}
&&\sum_{r=0}^{1-A_{ij}} 
\sum_{\eta=0}^{1-A_{ij}-r} 
\sum_{t=0}^r 
(-1)^r 
\biggl\lbrack {{1-A_{ij}}\atop{r}} \biggr\rbrack_i
\biggl\lbrack {{1-A_{ij}-r}\atop{\eta}} \biggr\rbrack_i
\biggl\lbrack {r\atop{t}} \biggr\rbrack_i
q_i^{(A_{ij}+\eta+2t)(1-A_{ij}-\eta-r)}
\\
&&
\qquad \qquad \qquad \times \quad 
q_i^{t(r-t)}
(q_i-q_i^{-1})^{\eta+t}
F_i^\eta 
F_j
F_i^t 
K_i^{A_{ij}+\eta+t-1}.
\end{eqnarray*}}\hfill
\parbox{1cm}{\begin{eqnarray}
\label{eq:c2}
\end{eqnarray}}

\noindent In 
(\ref{eq:c2})
we change variables by substituting
$u=r-t$ and $v=\eta+t$. We find
that for $0 \leq v \leq 1-A_{ij}$
and for $0 \leq t \leq v$
the coefficient of
$F_i^{v-t}F_jF_i^tK_i^{A_{ij}+v-1}$
in 
(\ref{eq:c2}) is equal to
\begin{eqnarray}
(-1)^t\biggl\lbrack {{1-A_{ij}}\atop{v}} \biggr\rbrack_i
\biggl\lbrack {{v}\atop{t}} \biggr\rbrack_i
q_i^{(A_{ij}+v+t)(1-A_{ij}-v)}
(q_i-q^{-1}_i)^v
\label{eq:c3}
\end{eqnarray}
times
\begin{eqnarray}
\sum_{u=0}^{1-A_{ij}-v} 
(-1)^{u} 
\biggl\lbrack {{1-A_{ij}-v}\atop{u}} \biggr\rbrack_i
q_i^{-u(A_{ij}+v)}.
\label{eq:c4}
\end{eqnarray}
The sum
(\ref{eq:c4}) is equal to 
$\delta_{v,1-A_{ij}}$
in view of
 (\ref{eq:vand2}).
For $v=1-A_{ij}$
the expression
(\ref{eq:c3}) is equal to
\begin{eqnarray}
(-1)^t 
\biggl\lbrack {{1-A_{ij}}\atop{t}} \biggr\rbrack_i
(q_i-q^{-1}_i)^{1-A_{ij}}.
\label{eq:c5}
\end{eqnarray}
Therefore 
(\ref{eq:c2}) is equal to
\begin{eqnarray}
(q_i-q^{-1}_i)^{1-A_{ij}}
\sum_{t=0}^{1-A_{ij}}
(-1)^t 
\biggl\lbrack {{1-A_{ij}}\atop{t}} \biggr\rbrack_i
F_i^{1-A_{ij}-t}F_jF_i^t
\label{eq:c6}
\end{eqnarray}
and this is zero in view of  
(R7). From these comments we find
(\ref{eq:b3}) is equal to zero.
We have now verified
(\ref{eq:hat7}).
Line 
(\ref{eq:hat8}) 
will follow from
(\ref{eq:hat7}) if we can display
an antiautomorphism
of 
$U_q(\mathfrak{g})$
that fixes
$
{\hat X}_i
$
and sends 
$
{\hat Y}_i
$
to 
$
{\hat Z}_i
$
for
$1 \leq i \leq n$.
The desired antiautomorphism will be the composition
of the following two maps.
Observe that there exists an antiautomorphism
of 
$U_q(\mathfrak{g})$ that satisfies
\begin{eqnarray}
K_i \rightarrow K_i,
\qquad \qquad 
E_i \rightarrow F_i,
\qquad \qquad 
F_i \rightarrow E_i
\label{eq:antione}
\end{eqnarray}
for $1 \leq i \leq n$.
Next observe that there exists an automorphism of
$U_q(\mathfrak{g})$ that satisfies
\begin{eqnarray}
K_i \rightarrow K_i,
\qquad \qquad
E_i \rightarrow -q_iK_i^{-1}E_i,
\qquad \qquad
F_i \rightarrow -q_i^{-1}F_iK_i
\label{eq:autone}
\end{eqnarray}
for $1 \leq i \leq n$.
Let $\mu $ denote the composition
of the map
(\ref{eq:antione}) followed by the map
(\ref{eq:autone}).
We observe $\mu$ is an antiautomorphism
of 
$U_q(\mathfrak{g})$. We routinely verify
that $\mu$ 
 fixes
$
{\hat X}_i
$
and sends 
$
{\hat Y}_i
$
to 
$
{\hat Z}_i
$
for
$1 \leq i \leq n$.
Applying $\mu$ to
the equation in
(\ref{eq:hat7}) we obtain
$(-1)^{1-A_{ij}}$ times the equation
in 
(\ref{eq:hat8}).
This shows that 
(\ref{eq:hat8}) holds in
$U_q(\mathfrak{g})$.
We have now shown that
 (\ref{eq:hat1})--(\ref{eq:hat8})
hold in
$U_q(\mathfrak{g})$.
Therefore 
(\ref{eq:i1})--(\ref{eq:i3}) gives
a homomorphism of 
$\K(q)$-algebras
from
$U$ to
$U_q(\mathfrak{g})$.
Next we show that
(\ref{eq:iinv1})--(\ref{eq:iinv3}) 
gives a homomorphism 
of $\K(q)$-algebras 
 from 
$U_q(\mathfrak{g})$ to
$U$.
In order to do this, we 
let ${\tilde E}_i, {\tilde F}_i, {\tilde K}^{\pm 1}_i$
denote the expressions on the right
in 
(\ref{eq:iinv1})--(\ref{eq:iinv3}) respectively
and show
the following hold in  
$U$:
\begin{eqnarray}
&&{\tilde K}_i{\tilde K}^{-1}_i = 
{\tilde K}^{-1}_i {\tilde K}_i =  1,
\label{eq:t1}
\\
&&{\tilde K}_i{\tilde K}_j = {\tilde K}_j {\tilde K}_i,
\label{eq:t2}
\\
&&{\tilde K}_i {\tilde E}_j{\tilde K}^{-1}_i=q_i^{A_{ij}}{\tilde E}_j,
\label{eq:t3}
\\
&&{\tilde K}_i{\tilde F}_j{\tilde K}^{-1}_i=
q_i^{-A_{ij}}{\tilde F}_j,
\label{eq:t4}
\\
&&{\tilde E}_i{\tilde F}_j-{\tilde F}_j{\tilde E}_i 
= \delta_{ij}\frac{{\tilde K}_i-{\tilde K}_i^{-1}}{q_i-q^{-1}_i},
\label{eq:t5}
\\
&&\sum_{r=0}^{1-A_{ij}} (-1)^r 
\biggl\lbrack {{1-A_{ij}}\atop{r}} \biggr\rbrack_i
{\tilde E}^{1-A_{ij}-r}_i {\tilde E}_j {\tilde E}^r_i =0 \qquad \mbox{if} \qquad i\not=j,
\label{eq:t6}
\\
&&
\sum_{r=0}^{1-A_{ij}} (-1)^r 
\biggl\lbrack {{1-A_{ij}}\atop{r}} \biggr\rbrack_i
{\tilde F}^{1-A_{ij}-r}_i {\tilde F}_j {\tilde F}^r_i =0 \qquad \mbox{if} \qquad i\not=j.
\label{eq:t7}
\end{eqnarray}
Line
(\ref{eq:t1}) holds by (E1)
and line 
(\ref{eq:t2}) holds by (E2).
Line 
(\ref{eq:t3}) follows
from (E1), (E2), (E4). 
Line 
(\ref{eq:t4}) follows
from (E1), (E2), (E3). 
Line 
(\ref{eq:t5}) follows
from (E1)--(E6).
Before we verify 
(\ref{eq:t6}) we make some comments.
In the equation of
(E4) we multiply each term on the left
by $X_j$ and simplify the result
using (E1), (E2) to obtain
\begin{eqnarray}
X_i(1-X_jZ_j)=q^{A_{ij}}_i(1-X_jZ_j)X_i.
%\qquad \qquad (1 \leq i,j\leq n).
\label{eq:xzcom}
\end{eqnarray}
In 
(\ref{eq:xzcom}) we set $i=j$ and $A_{ii}=2$ to obtain
\begin{eqnarray}
X_i(1-X_iZ_i)=q^2_i(1-X_iZ_i)X_i.
%\qquad \qquad (1 \leq i\leq n).
\label{eq:xzcom2}
\end{eqnarray}
Using
(\ref{eq:3ident}),
(\ref{eq:xzcom2}) and induction 
we obtain
\begin{eqnarray}
(1-X_iZ_i)^m
&=&
\sum_{r=0}^m (-1)^r
\biggl\lbrack {{m}\atop{r}} \biggr\rbrack_i
q_i^{r(1-m)} X_i^r Z_i^r 
\label{eq:xzp}
\end{eqnarray}
for all integers $m\geq 0$.
The expression on the left in
(\ref{eq:t6}) is a scalar multiple of
\begin{eqnarray}
\sum_{r=0}^{1-A_{ij}} (-1)^r 
\biggl\lbrack {{1-A_{ij}}\atop{r}} \biggr\rbrack_i
(1-X_iZ_i)^{1-A_{ij}-r} (1-X_jZ_j) (1-X_iZ_i)^r.
\label{eq:t6a}
\end{eqnarray}
We assume $i\not=j$ and show that
(\ref{eq:t6a}) is zero.
In (\ref{eq:t6a}) we evaluate
the terms $(1-X_iZ_i)^r$ 
using
(\ref{eq:xzp})
and obtain
\begin{eqnarray}
\sum_{r=0}^{1-A_{ij}}
(-1)^r 
\biggl\lbrack {{1-A_{ij}}\atop{r}} \biggr\rbrack_i
(1-X_iZ_i)^{1-A_{ij}-r} (1-X_jZ_j)
\sum_{\eta=0}^r (-1)^\eta
\biggl\lbrack {{r}\atop{\eta}} \biggr\rbrack_i
q_i^{\eta(1-r)} X_i^\eta Z_i^\eta.
\label{eq:xzp2}
\end{eqnarray}
In (\ref{eq:xzp2}) 
we evaluate the terms
$(1-X_jZ_j)X_i^\eta$ 
using
(\ref{eq:xzcom})
and obtain
\begin{eqnarray}
\sum_{r=0}^{1-A_{ij}}
\sum_{\eta=0}^r
(-1)^{r+\eta} 
\biggl\lbrack {{1-A_{ij}}\atop{r}} \biggr\rbrack_i
\biggl\lbrack {{r}\atop{\eta}} \biggr\rbrack_i
q_i^{\eta(1-r-A_{ij})}
(1-X_iZ_i)^{1-A_{ij}-r}
X_i^\eta(1-X_jZ_j)
Z_i^\eta.
\label{eq:xzp4}
\end{eqnarray}
In  
(\ref{eq:xzp4}) we evaluate
the terms
$(1-X_iZ_i)^{1-A_{ij}-r}
X_i^\eta$
using
(\ref{eq:xzcom2}) and obtain
\begin{eqnarray}
\sum_{r=0}^{1-A_{ij}}
\sum_{\eta=0}^r
(-1)^{r+\eta} 
\biggl\lbrack {{1-A_{ij}}\atop{r}} \biggr\rbrack_i
\biggl\lbrack {{r}\atop{\eta}} \biggr\rbrack_i
q_i^{\eta(A_{ij}+r-1)}
X_i^\eta
(1-X_iZ_i)^{1-A_{ij}-r}
(1-X_jZ_j)
Z_i^\eta.
\label{eq:xzp5}
\end{eqnarray}
Observe 
(\ref{eq:xzp5}) is equal to
\begin{eqnarray}
&&\sum_{r=0}^{1-A_{ij}}
\sum_{\eta=0}^r
(-1)^{r+\eta} 
\biggl\lbrack {{1-A_{ij}}\atop{r}} \biggr\rbrack_i
\biggl\lbrack {{r}\atop{\eta}} \biggr\rbrack_i
q_i^{\eta(A_{ij}+r-1)}
X_i^\eta
(1-X_iZ_i)^{1-A_{ij}-r}
Z_i^\eta
\label{eq:xzp6}
\end{eqnarray}
minus
\begin{eqnarray}
\sum_{r=0}^{1-A_{ij}}
\sum_{\eta=0}^r
(-1)^{r+\eta} 
\biggl\lbrack {{1-A_{ij}}\atop{r}} \biggr\rbrack_i
\biggl\lbrack {{r}\atop{\eta}} \biggr\rbrack_i
q_i^{\eta(A_{ij}+r-1)}
X_i^\eta
(1-X_iZ_i)^{1-A_{ij}-r}
X_jZ_j
Z_i^\eta.
\label{eq:xzp7}
\end{eqnarray}
So far we have shown
that 
(\ref{eq:t6a}) is equal to
(\ref{eq:xzp6}) minus 
(\ref{eq:xzp7}).
We  evaluate 
(\ref{eq:xzp6}), 
(\ref{eq:xzp7}) separately.
In 
(\ref{eq:xzp6}) we eliminate
the terms $(1-X_iZ_i)^{1-A_{ij}-r}$
using
(\ref{eq:xzp})
and obtain
\begin{eqnarray}
\sum_{r=0}^{1-A_{ij}}
\sum_{\eta=0}^r
\sum_{t=0}^{1-A_{ij}-r}
(-1)^{r+\eta+t} 
\biggl\lbrack {{1-A_{ij}}\atop{r}} \biggr\rbrack_i
\biggl\lbrack {{r}\atop{\eta}} \biggr\rbrack_i
\biggl\lbrack {{1-A_{ij}-r}\atop{t}} \biggr\rbrack_i
q_i^{\eta(A_{ij}+r-1)}
q_i^{t(A_{ij}+r)}
X_i^{\eta+t}
Z_i^{\eta+t}.
\label{eq:xzp8}
\end{eqnarray}
In 
(\ref{eq:xzp8}) we change variables
by substituting $r=u+\eta$ and $t=v-\eta$.
We find that 
for $0 \leq v \leq 1-A_{ij}$ the
coefficient of
$X_i^vZ_i^v$ in 
(\ref{eq:xzp8}) is equal to
\begin{eqnarray}
(-1)^v
\biggl\lbrack {{1-A_{ij}}\atop{v}} \biggr\rbrack_i
q_i^{v A_{ij}}
\label{eq:xzp9}
\end{eqnarray}
times
\begin{eqnarray}
\sum_{\eta=0}^v 
(-1)^\eta 
\biggl\lbrack {{v}\atop{\eta}} \biggr\rbrack_i
q_i^{\eta(v-1)}
\label{eq:xzp10}
\end{eqnarray}
times
\begin{eqnarray}
\sum_{u=0}^{1-A_{ij}-v} 
(-1)^{u} 
\biggl\lbrack {{1-A_{ij}-v}\atop{u}} \biggr\rbrack_i
q_i^{uv}.
\label{eq:xzp11}
\end{eqnarray}
The sum (\ref{eq:xzp10}) is equal to
$\delta_{0,v}$ in view of
 (\ref{eq:vand2}).
For $v=0$ the expression
(\ref{eq:xzp9}) is equal to 1 and
the sum (\ref{eq:xzp11}) is equal to
\begin{eqnarray}
\label{eq:xzp12}
 \prod_{s=0}^{-A_{ij}}(1-q_i^{A_{ij}+2s})
\end{eqnarray}
in view of
(\ref{eq:vand2}).
Therefore
(\ref{eq:xzp6}) is equal to
(\ref{eq:xzp12}).
In 
(\ref{eq:xzp7}) we evaluate
the terms
$(1-X_iZ_i)^{1-A_{ij}-r}X_j$
using
(\ref{eq:xzcom}) and simplify the result using
$q_i^{A_{ij}}=q_j^{A_{ji}}$.
We obtain
\begin{eqnarray}
\sum_{r=0}^{1-A_{ij}}
\sum_{\eta=0}^r
(-1)^{r+\eta} 
\biggl\lbrack {{1-A_{ij}}\atop{r}} \biggr\rbrack_i
\biggl\lbrack {{r}\atop{\eta}} \biggr\rbrack_i
q_i^{(A_{ij}+\eta)(A_{ij}+r-1)}
X_i^\eta
X_j
(1-X_iZ_i)^{1-A_{ij}-r}
Z_j
Z_i^\eta.
\label{eq:xzp13}
\end{eqnarray}
In 
(\ref{eq:xzp13}) we evaluate the
terms
$(1-X_iZ_i)^{1-A_{ij}-r}$ using
(\ref{eq:xzp}) and obtain

\parbox{10cm}
{\begin{eqnarray*}
&&\sum_{r=0}^{1-A_{ij}}
\sum_{\eta=0}^r
\sum_{t=0}^{1-A_{ij}-r}
(-1)^{r+\eta+t} 
\biggl\lbrack {{1-A_{ij}}\atop{r}} \biggr\rbrack_i
\biggl\lbrack {{r}\atop{\eta}} \biggr\rbrack_i
\biggl\lbrack {{1-A_{ij}-r}\atop{t}} \biggr\rbrack_i
\\
&& \qquad \qquad \qquad \qquad \qquad
\times \quad 
q_i^{(A_{ij}+\eta)(A_{ij}+r-1)}
q_i^{t(A_{ij}+r)}
X_i^{\eta+t}
X_j
Z_i^t
Z_j
Z_i^\eta.
\end{eqnarray*}}\hfill
\parbox{1cm}
{\begin{eqnarray}
\label{eq:xzp14}
\end{eqnarray}}

\noindent
In
(\ref{eq:xzp14}) 
we change variables by substituting
$r=u+\eta$ and $t=v-\eta$. We find
that for $0 \leq \eta \leq 1-A_{ij}$
and for $\eta \leq v \leq 1-A_{ij}$
the coefficient of 
$X_i^vX_j
Z_i^{v-\eta}Z_jZ_i^\eta$ in 
(\ref{eq:xzp14}) is equal to
\begin{eqnarray}
(-1)^{v+\eta}
\biggl\lbrack {{1-A_{ij}}\atop{v}} \biggr\rbrack_i
\biggl\lbrack {{v}\atop{\eta}} \biggr\rbrack_i
q_i^{(A_{ij}+\eta)(A_{ij}+v-1)}
%X_i^vX_j
\label{eq:xzp15}
\end{eqnarray}
times
\begin{eqnarray}
\sum_{u=0}^{1-A_{ij}-v}
(-1)^{u} 
\biggl\lbrack {{1-A_{ij}-v}\atop{u}} \biggr\rbrack_i
q_i^{u(A_{ij}+v)}.
\label{eq:xzp16}
\end{eqnarray}
The sum
(\ref{eq:xzp16}) 
is equal to $\delta_{1-A_{ij},v}$
in view of
 (\ref{eq:vand2}).
For 
$v=1-A_{ij}$
the expression
(\ref{eq:xzp15}) is equal to
\begin{eqnarray}
(-1)^{1-A_{ij}+\eta}
\biggl\lbrack {{1-A_{ij}}\atop{\eta}} \biggr\rbrack_i.
%X_i^{1-A_{ij}}X_j.
\label{eq:xzp17}
\end{eqnarray}
From these comments we find
(\ref{eq:xzp14}) is equal to
\begin{eqnarray}
(-1)^{1-A_{ij}}
X_i^{1-A_{ij}}X_j
\sum_{\eta=0}^{1-A_{ij}}
(-1)^\eta
\biggl\lbrack {{1-A_{ij}}\atop{\eta}} \biggr\rbrack_i
Z_i^{1-A_{ij}-\eta}Z_jZ_i^\eta.
\label{eq:xzp18}
\end{eqnarray}
Using
(E8) we find
(\ref{eq:xzp18}) is equal to
\begin{eqnarray}
\label{eq:xzp19}
(-1)^{1-A_{ij}}\prod_{s=0}^{-A_{ij}}(1-q_i^{A_{ij}+2s}).
\end{eqnarray}
In line  
(\ref{eq:xzp19}) 
the product
is 0 if $A_{ij}$ is even
and of course
$(-1)^{1-A_{ij}}=1$ if $A_{ij}$ is odd.
Therefore
(\ref{eq:xzp19}) 
is equal to
(\ref{eq:xzp12}).
We have now shown that
(\ref{eq:xzp7}) is equal to
(\ref{eq:xzp12}).
Since each of
(\ref{eq:xzp6}),
(\ref{eq:xzp7}) is equal 
to 
(\ref{eq:xzp12}), their difference is zero.
We showed earlier that
(\ref{eq:t6a}) is equal to this difference
so
(\ref{eq:t6a}) is equal to zero.
We have now verified (\ref{eq:t6}).
Before we verify 
(\ref{eq:t7}) we have a comment.
Observe that there exists an antiautomorphism
of $U$ that satisfies
\begin{eqnarray}
X_i \rightarrow X_i,
\qquad \qquad 
Y_i \rightarrow Z_i,
\qquad \qquad 
Z_i \rightarrow Y_i
\label{eq:backanti}
\end{eqnarray}
for $1 \leq i \leq n$.
The expression on the left in
(\ref{eq:t7}) is a scalar multiple of
\begin{eqnarray}
\sum_{r=0}^{1-A_{ij}} (-1)^r 
\biggl\lbrack {{1-A_{ij}}\atop{r}} \biggr\rbrack_i
(Y_i-X^{-1}_i)^{1-A_{ij}-r} (Y_j-X^{-1}_j) (Y_i-X^{-1}_i)^r.
\label{eq:backt7a}
\end{eqnarray}
We assume $i\not=j$ and show that
(\ref{eq:backt7a}) is zero.
Earlier we showed that
(\ref{eq:t6a}) is zero.
Applying the antiautomorphism
(\ref{eq:backanti})  to 
(\ref{eq:t6a}) we find
\begin{eqnarray}
\sum_{r=0}^{1-A_{ij}} (-1)^r 
\biggl\lbrack {{1-A_{ij}}\atop{r}} \biggr\rbrack_i
(1-Y_iX_i)^r (1-Y_jX_j) (1-Y_iX_i)^{1-A_{ij}-r} =0.
\label{eq:backt8}
\end{eqnarray}
In
(\ref{eq:backt8}) we change variables by replacing
$r$ by $1-A_{ij}-r$. We also multiply both sides
by $-1$. We find 
\begin{eqnarray}
\sum_{r=0}^{1-A_{ij}} (-1)^r 
\biggl\lbrack {{1-A_{ij}}\atop{r}} \biggr\rbrack_i
(Y_iX_i-1)^{1-A_{ij}-r}(Y_jX_j-1) (Y_iX_i-1)^r=0.
\label{eq:backt9}
\end{eqnarray}
Applying the antiautomorphism
(\ref{eq:backanti}) 
to 
(\ref{eq:xzcom}), 
(\ref{eq:xzcom2}) and manipulating
the result we obtain
\begin{eqnarray}
(Y_jX_j-1)X^{-1}_i
&=&
q^{-A_{ij}}_iX^{-1}_i(Y_jX_j-1),
\label{eq:backxzcom}
\\
(Y_iX_i-1)X^{-1}_i&=&q^{-2}_iX^{-1}_i(Y_iX_i-1).
\label{eq:backxzcom2}
\end{eqnarray}
Multiplying both sides of
(\ref{eq:backt9}) on the right
by $X_i^{A_{ij}-1}X_j^{-1}$
and evaluating the result using
(\ref{eq:backxzcom}), 
(\ref{eq:backxzcom2})
we routinely find that
(\ref{eq:backt7a}) is zero.
By the above comments
(\ref{eq:t7}) holds in
$U$. We have now shown that
(\ref{eq:t1})--(\ref{eq:t7})  hold 
in $U$.
Therefore 
(\ref{eq:iinv1})--(\ref{eq:iinv3}) 
gives a homomorphism of $\K(q)$-algebras
 from 
$U_q(\mathfrak{g})$ to
$U$.
So far we have shown
that 
(\ref{eq:i1})--(\ref{eq:i3}) 
gives a homomorphism of $\K(q)$-algebras
 from 
$U$ to
$U_q(\mathfrak{g})$ 
and that
(\ref{eq:iinv1})--(\ref{eq:iinv3}) 
gives a homomorphism of $\K(q)$-algebras
 from 
$U_q(\mathfrak{g})$ 
to $U$.
One routinely verifies that these maps 
are inverses. Therefore, each of these maps
is a bijection and hence an isomorphism
of $\K(q)$-algebras.
\hfill $\Box $ \\

\begin{definition}
\rm
By the {\it equitable presentation} for
$U_q(\mathfrak{g})$ we mean
the presentation given in
 Theorem
\ref{thm:uqg2}.
We call the generators
$X_i, X^{-1}_i, Y_i, Z_i$ $(i =1,\ldots, n)$
the {\it equitable generators}.
\end{definition}

\begin{definition}
\label{def:ncg}
\rm
For notational convenience, we identify
the copy of 
$U_q(\mathfrak{g})$
given 
in Definition
\ref{def:uqg} with the copy of
$U_q(\mathfrak{g})$
given in Theorem
\ref{thm:uqg2}, via the isomorphism
given 
in Theorem
\ref{thm:uqg2}.
\end{definition}

\noindent
The Hopf algebra 
structure for 
$U_q(\mathfrak{g})$
given in Lemma
\ref{lem:hopfuqg}
looks as follows 
in terms of the 
equitable generators.

\begin{theorem}
With reference to
 Lemma
\ref{lem:hopfuqg} and
Definition \ref{def:ncg},
the comultiplication $\Delta$ 
 satisfies
\begin{eqnarray*}
\Delta(X_i) &=& X_i \otimes X_i,
\\
\Delta(Y_i) &=& 
(Y_i-1) \otimes X_i^{-1}+
1 \otimes Y_i,
\\
\Delta(Z_i) &=& 
(Z_i-1) \otimes X_i^{-1}+
1 \otimes Z_i.
\end{eqnarray*}
The counit $\varepsilon$ 
satisfies
\begin{eqnarray*}
\varepsilon(X_i)=1, \qquad \qquad \varepsilon(Y_i)=1,
\qquad \qquad 
\varepsilon(Z_i)=1.
\end{eqnarray*}
The antipode $S$ 
 satisfies
\begin{eqnarray*}
S(X_i)= X_i^{-1},
\qquad \qquad 
S(Y_i)= 1+X_i-Y_iX_i,
\qquad \qquad 
S(Z_i) = 1+X_i-Z_iX_i.
\end{eqnarray*}
\end{theorem}
\noindent {\it Proof:} 
Routine verification.
\hfill $\Box $ \\

\noindent We finish this section with some remarks.

\begin{remark}
\rm
Referring to the relations (E7), (E8)
in Theorem 
\ref{thm:uqg2}, if $A_{ij}$ is
even then the product on the right-hand side is
zero.
\end{remark}

\begin{remark}
\label{koorn3}
\rm
With reference to
Lemma
\ref{lem:hopfuqg}
and Definition \ref{def:ncg},
for $1 \leq i \leq n$ and for
$y_i=Y_i-1$, $z_i=Z_i-1$ we have
\begin{eqnarray*}
\Delta(y_i) &=& 
y_i\otimes X_i^{-1}+
1 \otimes y_i,
\\
\Delta(z_i) &=& 
z_i\otimes X_i^{-1}+
1 \otimes z_i.
\end{eqnarray*}
We will discuss these two equations in
the next section.
%Therefore, each of 
%$y_i, z_i$ is a twisted primitive element
%with respect to $X_i$, in the sense of
%Koornwinder
%\cite[p.~801]{koorn}.
\end{remark}

\begin{remark}
\rm
For the quantum group of type $A^{(1)}_1$
the equitable presentation is essentially
the same as the presentation
given in
\cite[Theorem 2.1]{tdanduq}.
\end{remark}

\section{The equitable generators and twisted primitive elements}

In this section we discuss how 
the equitable generators are related to Koornwinder's
twisted primitive elements
\cite{koorn1},
\cite{koorn}.
We begin with a definition.

\begin{definition}
\label{def:glike}
\rm
\cite[p. 56]{Kassel}
 For a Hopf algebra $U$, an element $u \in U$
is called {\it group-like} whenever $u \not=0$ and
$\Delta(u)=u\otimes u$.
We let $G(U)$ denote the set of group-like elements
of $ U$.
\end{definition}

\begin{lemma}
\cite[p. 56]{Kassel}
 With reference to
Definition
\ref{def:glike},
for all $u,v \in 
 G(U)$ we have
$uv \in G(U)$.
  Also for $u \in 
 G(U)$, $u$ has an inverse in 
 $G(U)$ which is equal to $S(u)$.
Consequently 
 $G(U)$ is a group.
\end{lemma}

\begin{example}
\rm
\cite[Lemma 1]{masuda}
For
the Hopf algebra
$U=U_q(\mathfrak{sl}_2)$ the group
$G(U)$ consists of the elements
$X^i$  $(i \in \Z)$, where $X$ is
from Definition
\ref{def:equit2}.
\end{example}

\begin{definition}
\rm
\cite[p. 48]{Kassel}
 With reference to
Definition
\ref{def:glike}, 
an element $u \in 
U$ is called {\it primitive}
whenever 
$\Delta(u)=
1\otimes u+
u\otimes 1$.
\end{definition}

\begin{definition}
\rm
\cite[p.~801]{koorn}
With reference to
Definition
\ref{def:glike}, 
for $u \in 
U$ and for 
$g \in G(U)$ we say
$u$ is {\it twisted primitive with respect to $g$}
whenever 
$\Delta(u)=g\otimes u+u\otimes S(g)$.
\end{definition}

\begin{example}
\label{ex:uqsl2}
\rm
For the Hopf algebra
$U_q(\mathfrak{sl}_2)$
an element $u$ is twisted primitive with respect to
$X$ whenever 
$\Delta(u)=
X\otimes u+
u\otimes X^{-1}$.
\end{example}

\noindent
Comparing Example
\ref{ex:uqsl2} with
 Remark
\ref{koorn}
we find that although $y,z$ are not
twisted primitive with respect to
$X$, there is some resemblence.
In order to interpret this resemblence we make a definition.

\begin{definition}
\label{def:qtp}
\rm With reference to
Definition
\ref{def:glike}, 
for $u \in U$ and for $g \in G(U)$
we say $u$ is {\it quasi twisted primitive with respect
to $g$} whenever 
$\Delta(u) = 1 \otimes u + u \otimes S(g)$.
\end{definition}

\noindent
The following lemma shows how
quasi twisted primitive elements are related to
twisted primitive elements. The proof is routine and
omitted.

\begin{lemma}
\label{lem:bij}
With reference to
Definition
\ref{def:glike}, 
for $g \in G(U)$ the map
$u \mapsto uS(g)$ gives a bijection
from (i) the set of elements in
$U$ that are twisted primitive with respect to $g$,
to (ii) the set of elements in $U$ that
are quasi twisted primitive with respect to $g^2$.
\end{lemma}

\begin{remark}
\rm With reference to
Definition \ref{def:glike}, 
in general the set of squares 
$\lbrace g^2 \,|\,g \in G(U)\rbrace$
is a proper subset of $G(U)$. By this and
in view of Lemma
\ref{lem:bij}, 
the concept of
a quasi twisted primitive element is a bit more general
than the concept of a twisted primitive element.
\end{remark}

\begin{theorem}
Referring to Remark
\ref{koorn},
each of $y,z$
is quasi twisted primitive with respect to $X$.
Referring to Remark
\ref{koorn3}, for $1 \leq i \leq n$
each of 
$y_i, z_i$ is quasi twisted primitive with respect
to $X_i$.
\end{theorem}
\noindent {\it Proof:} 
Immediate from
 Remark
\ref{koorn}, 
 Remark
\ref{koorn3}, and
Definition
\ref{def:qtp}.
\hfill $\Box $ \\

\noindent
We refer the reader to
\cite{koe1},
\cite{koe2},
\cite{koorn1},
\cite{noumi4},
\cite{nd},
\cite{noumi3},
\cite{noumi},
\cite{ros},
\cite{stokman} for more information about
twisted primitive elements and related topics.

%%%%%%%%

\bigskip
\noindent Paul Terwilliger \hfil\break
Department of Mathematics \hfil\break
University of Wisconsin \hfil\break
480 Lincoln Drive \hfil\break
Madison, Wisconsin, 53706-1388 USA 
\hfil\break
{Email: \tt terwilli@math.wisc.edu}

\end{document}